\documentclass{svmult}
\usepackage{amsmath,amssymb,bbm}
\smartqed

\begin{document}

\title*{Convergence rates for density estimators of weakly dependent time series}
\author{Nicolas Ragache\inst{1}\and
Olivier Wintenberger\inst{2}}
\institute{MAP5, Universit\'e Ren\'e Descartes 45 rue des
Saints-P\`{e}res, 75270 Paris, France \texttt{nicolas.ragache@ensae.fr} \and
SAMOS, Statistique Appliqu\'ee et MOd\'elisation Stochastique,
Universit\'e Paris 1, Centre Pierre Mend\`es France, 90 rue de
Tolbiac, F-75634 Paris Cedex 13, France.
\texttt{olivier.wintenberger@univ-paris1.fr}}

\maketitle

\section{Introduction}
\label{sec:1}
Assume that $(X_n)_{n\in\mathbb{Z}}$ is a sequence of $\mathbb{R}^d$
valued random variables with common distribution which is absolutely
continuous with respect to Lebesgue's measure, with density $f$.
Stationarity is not assumed so that the case of a sampled process
$\{X_{i,n}=x_{h_n(i)}\}_{1\le i\le n}$ for any sequence of monotonic
functions $(h_n(.))_{n\in\mathbb{Z}}$ and any stationary process
$(x_n)_{n\in\mathbb{Z}}$ that admits a marginal density is included.
This paper investigates convergence rates for density estimation in
different cases.  First, we consider two concepts of weak dependence:
\begin{itemize}
\item Non-causal $\eta$-dependence introduced in \cite{dou2} by
  Doukhan \& Louhichi,
\item Dedecker \& Prieur's $\tilde \phi$-dependence (see \cite{ded2}).
\end{itemize}
These two notions of dependence cover a large number of examples of
time series (see section \S~3). Next, following Doukhan (see
\cite{dou1}) we propose a unified study of linear density estimators
$\hat f_n$ of the form
\begin{equation}\label{estlin}
\hat f_n(x) = \frac1{n} \sum_{i=1}^n K_{m_n}(x,X_i)\;,
\end{equation}
where $\{K_{m_n}\}$ is a sequence of kernels. Under classical
assumptions on $\{K_{m_n}\}$ (see section \S~2.2), the results in
the case of independent and identically distributed (i.i.d. in
short) observations $X_i$ are well known (see for instance
\cite{t}). At a fixed point $x\in\mathbb{R}^d$, the sequence $m_n$
can be chosen such that
\begin{equation}\label{mini1}
\|\hat f_n(x)-f(x)\|_q = O \left(n^{-\rho/(2\rho+d)}\right)\;,
\end{equation}

\noindent where $\|X\|_q^q=\mathbb{E}|X|^q$. The coefficient
$\rho>0$ measures the regularity of $f$ (see Section~\ref{nes} for
the definition of the notion of regularity). The same rate of
convergence also holds for the Mean Integrated Square Error (MISE),
defined as $\int\| \hat f_n(x) - f(x)\|_2^2 p(x) \, \D x$ for some
nonnegative and integrable function $p$.  The rate of uniform
convergence on a compact set incurs a logarithmic loss appears. For
all $M>0$ and for a suitable choice of the sequence $m_n$,
\begin{equation}\label{mini2}
\mathbb{E}\sup_{\|x\|\le M}|\hat f_n(x)-f(x)|^q={\cal O}\left(\frac{\log
n}n\right)^{q\rho/(d+2\rho)}\;,
\end{equation}
and
\begin{equation}\label{mini3}
\sup_{\|x\|\le M}|\hat f_n(x)-f(x)|=_{a.s.}{\cal O}\left(\frac{\log
n}n\right)^{\rho/(d+2\rho)}\;.
\end{equation}
These rates are optimal in the minimax sense. We thus have no hope to
improve on them in the dependent setting.  A wide literature deals
with density estimation for absolutely regular or $\beta$-mixing
processes (for a definition of mixing coefficients, see \cite{dou0}).
For instance, under the assumption $\displaystyle \beta_r=o
\left(r^{-3-2d/\rho}\right)$, Ango Nze \& Doukhan prove in \cite{ang0}
that (\ref{mini1}), (\ref{mini2}) and (\ref{mini3}) still hold. The
sharper condition $ \sum_r|\beta_r|<\infty$ entails the optimal rate
of convergence for the MISE (see \cite{v}). Results for the MISE have
been extended to the more general $\tilde \phi$- and $\eta$-dependence
contexts by Dedecker \& Prieur (\cite{ded2}) and Doukhan \& Louhichi
in \cite{dou3}.  In this paper, our aim is to extend the bounds
(\ref{mini1}), (\ref{mini2}) and (\ref{mini3}) in the $\eta$- and
$\tilde \phi$-weak dependence contexts.

 We use the same method as in \cite{dou2} based on the
following moment inequality for weakly dependent and centered
sequences $(Z_n)_{n\in\mathbb{Z}}$. For each even integer $q$ and for each
integer $n\ge 2$:
\begin{equation}\label{douklou}
\left\|\sum_{i=1}^nZ_i\right\|_q^q\le
\frac{(2q-2)!}{(q-1)!}\left\{V_{2,n}^{q/2}\vee V_{q,n}\right\}\;,
\end{equation}
where $\|X\|_q^q=\mathbb{E}|X|^q$ and for $k = 2, \dots,q$,
$$
V_{k,n} = n \sum_{r=0}^{n-1}(r+1)^{k-2} C_k(r) \; ,
$$
with
\begin{equation} \label{ckr}
C_{k}(r):=\sup\{|\mathrm{cov}(Z_{t_1}\cdots
Z_{t_p},Z_{t_{p+1}}\cdots Z_{t_k})|\}\;,
\end{equation}
where the supremum is over all the ordered $k$-tuples $t_1\le\cdots\le
t_k$ such that $\sup_{1\le i\le k-1}t_{i+1}-t_i=r$.

We will apply this bound when the $Z_i$s are defined in such a way
that $\sum_{i=1}^nZ_i$ is proportional to the fluctuation term $\hat
f_n(x)-\mathbb{E} \hat f_n(x)$. The inequality (\ref{douklou}) gives a
bound for this part of the deviation of the estimator which depends on
the covariance bounds $C_{k}(r)$.  The other part of the deviation is
the bias, which is treated by deterministic methods. In order to
obtain suitable controls of the fluctuation term, we need two
different type of bounds for $C_{k}(r)$.  Conditions on the decay of
the weak dependence coefficients give a first bound. Another type of
condition is also required to bound $C_{k}(r)$ for the smaller values
of $r$; this is classically achieved with a regularity condition on
the joint law of the pairs $(X_j,X_k)$ for all $j \ne k$. In Doukhan
\& Louhichi (see \cite{dou3}), rates of convergence are obtained when
the coefficient $\eta$ decays geometrically fast and the joint
densities are bounded. We relax these conditions to cover the case
when the joint distributions are not absolutely continuous and when
the $\eta$- and $\tilde \phi$-dependence coefficients decrease slowly
(sub-geometric and Riemannian decays are considered).

Under our assumptions, we prove that (\ref{mini1}) still holds (see
Theorem~\ref{theorem:1}).  Unfortunately, additional losses appear for
the uniform bounds. When $\eta_r$ or $\tilde\phi_r= O (e^{-ar^b})$
with $a>0$ and $b>0$, we prove in Theorem~\ref{thorem:2} that
(\ref{mini2}) and (\ref{mini3}) hold with $\log(n)$ replaced by
$\log^{2(b+1)/b}(n)$.  If $\eta_r$ or $\tilde\phi_r =
O(r^{-a})$ with $a>1$, Theorem~\ref{theorem:3} gives bounds similar to
(\ref{mini2}) and (\ref{mini3}) with the right hand side replaced by
$O ( n^{-q\rho/\{d+2\rho+2d/(q_0+d)\}}$ and $O (
\{\log^{q_0+d}(n)/n^{q_0-2} \}^{\rho/\{2\rho q_0+d(q_0+2)\}} )$,
respectively, and with $q_0 = 2 \lceil (a-1)/2 \rceil$ (by definition
$\lceil x \rceil$ is the smallest integer larger than or equal to the
real number $x$). As already noticed in \cite{dou3}, the loss w.r.t
the i.i.d. case highly depends on the decay of the dependence
coefficients. In the case of geometric decay, the loss is logarithmic
while it is polynomial in the case of polynomial decays.

\bigskip The paper is organized as follows. In
Section~\ref{sec:weakdependence}, we introduce the notions of $\eta$
and $\tilde \phi$ dependence. We give the notation and hypothesis in
Section~\ref{nes}. The main results are presented in
Section~\ref{sec:results}. We then apply these results to particular
cases of weak dependence processes, and we provide examples of
kernel $K_m$ in Section~\ref{sec:applications}.
Section~\ref{sec:proof} contains the proof of the Theorems and three
important lemmas.

\section{Main results}
We first describe the notions of dependence considered in this paper,
then we introduce assumptions and formulate the main results of the
paper (convergence rates).

\subsection{Weak dependence} \label{sec:weakdependence}
We consider a sequence $(X_i)_{i\in\mathbb{Z}}$ of $\mathbb{R}^d$
valued random variables, and we fix a norm $\|\cdot\|$ on
$\mathbb{R}^d$.  Moreover, if $h:\mathbb{R}^{du}\to \mathbb{R}$ for
some $u\ge1$, we define
\begin{eqnarray*}
\mbox{Lip }(h) & = & \sup_{(a_1,\ldots,a_u)\ne
(b_1,\ldots,b_u)}\frac{\left|h(a_1,\ldots,a_u)- h
(b_1,\ldots,b_u)\right|}{\|a_1-b_1\|+\cdots+\|a_u-b_u\|}\;.
\end{eqnarray*}

\begin{definition}[$\eta$-dependence, Doukhan  \& Louhichi (1999)]
    The process $(X_i)_{i\in\mathbb{Z}}$ is $\eta$-weakly dependent if
  there exists a sequence of non-negative real numbers $(\eta_r)_{r\ge
    0}$ satisfying $\eta_r\to0$ when $r\to\infty$ and
\begin{eqnarray*}
\left|\mathrm{cov}\left(h\left(X_{i_1}, \ldots X_{i_u}\right),
k\left(X_{i_{u + 1}},\ldots, X_{i_{u + v}}\right)\right)\right|
&\leq& (u{\rm Lip}(h)+v{\rm Lip}(k)) \eta_r\;,
\end{eqnarray*} for all  $(u+v)$-tuples, $(i_1, \ldots
,i_{u+v})$ with $i_1 \leq \cdots \leq i_u \le i_u+r \leq i_{u+1}\le
\cdots \le i_{u+v}$, and $h,k \in \Lambda^{(1)}$ where
\begin{eqnarray*}
 {\Lambda}^{(1)} &= &\left\{ h:\exists u\ge0,
 h:\mathbb{R}^{du}\to\ \mathbb{R},  \mathrm{Lip}\, (h) < \infty,
\|h\|_\infty=\sup_{x\in\mathbb{R}^{du}}|h(x)|\le 1
\right\} \; .
\end{eqnarray*}
\end{definition}

\noindent{\bf Remark}
The $\eta$-dependence condition can be applied to non-causal sequences
because information ``from the future'' (i.e. on the right of the
covariance) contributes to the dependence coefficient in the same way
as information ``from the past'' (i.e.  on the left). It is the
non-causal alternative to the $\theta$ condition in \cite{ded1} and
\cite{dou2}.

\begin{definition}[$\tilde\phi$-dependence,  Dedecker \& Prieur (2004)]
  Let $(\Omega,\mathcal A,\mathbb{P})$ be a probability space and
  $\mathcal M$ a $\sigma$-algebra of $\mathcal A$. For any
  $l\in\mathbb{N}^*$, any random variable $X\in\mathbb{R}^{dl}$ we
  define:
  $$
  \tilde \phi(\mathcal M,X)=\sup\{\|\mathbb{E}(g(X)|\mathcal
  M)-\mathbb{E}(g(X))\|_\infty,g\in \Lambda_{1,l}\}\;,
  $$
  where
  $\Lambda_{1,l}=\{h:\mathbb{R}^{dl}\mapsto\mathbb{R}/\mathrm{Lip}\,(h)<1\}$.
  The sequence of coefficients $\tilde \phi_{k}(r)$ is then defined by
  $$
  \tilde \phi_{k}(r)=\max_{l\le k}\frac{1}{l}\sup_{i+r\le
    j_1<j_2<\cdots<j_l}\tilde \phi(\sigma(\{X_j; j\le
  i\}),(X_{j_1},\dots,X_{j_l}))\;.
  $$
  The process is $\tilde \phi$-dependent if $\tilde
  \phi(r)=\sup_{k>0}\tilde \phi_{k}(r)$ tends to $0$ with $r$.
\end{definition}

\noindent {\bf Remark} The $\tilde \phi$ dependence coefficients provide
covariance bounds. For a Lipschitz function $k$ and a bounded function
$h$,
\begin{multline}
  \left| \mbox{cov} \left(h \left(X_{i_1},\dots, X_{i_u}\right), k
      \left(X_{i_{u + 1}},\dots, X_{i_{u + v}}\right) \right) \right|
  \\
  \leq v\mathbb{E}\left|h\left(X_{i_1},\dots, X_{i_u}\right)\right|
  \mathrm{Lip}\,(k) \tilde \phi(r)\; . \label{tau}
\end{multline}

\subsection{Notations and definitions}
\label{nes}

\noindent Assume that $(X_n)_{n\in\mathbb{Z}}$ is an  $\eta$ or $\tilde
\phi$ dependent sequence of $\mathbb{R}^d$ valued random variables. We
consider two types of decays for the coefficients. The geometric case
is the case when  Assumption [H1] or [H1'] holds.
\begin{itemize}
\item[] [H1]: $\eta_r = O \left(e^{-ar^b}\right)$ with $a>0$ and
  $b>0$,
\item[] [H1']: $\tilde \phi(r)= O \left(e^{-ar^b}\right)$ with $a>0$
  and $b>0$.
\end{itemize}
The Riemannian case is the case when Assumption [H2] or [H2'] holds.
\begin{itemize}
\item[] [H2]: $\eta_r={O}(r^{-a})$ with
$a>1$,
\item[] [H2']: $\tilde \phi(r)={\cal
O}(r^{-a})$ with $a>1$.
\end{itemize}

\noindent
As usual in density estimation, we shall assume:
\begin{itemize}
\item[] [H3]: The common marginal distribution of the random variables
  $X_n$, $n\in \mathbb{Z}$ is absolutely continuous with respect to
  Lebesgue's measure, with common bounded density $f$.
\end{itemize}

\noindent
The next assumption is on the density with respect to Lebesgue's
measure (if it exists) of the joint distribution of the pairs
$(X_j,X_k)$, $j \ne k$.
\begin{itemize}
\item[] [H4] The density $f_{j,k}$ of the joint distribution of the
  pair $(X_j,X_k)$ is uniformly bounded with respect to $j \ne k$.
\end{itemize}
\noindent Unfortunately, for some processes, these densities may not
even exist. For example, the joint distributions of Markov chains
$X_n=G(X_{n-1},\epsilon_n)$ may not be absolutely continuous. One of
the simplest example is
\begin{equation}\label{markov}
X_k=\frac12\left(X_{k-1}+\epsilon_k\right)\;,
\end{equation}
where $\{\epsilon_k\}$ is an i.i.d. sequence of Bernoulli random
variables and $X_0$ is uniformly distributed on $[0,1]$.  The process
$\{X_n\}$ is strictly stationary but the joint distributions of the
pairs $(X_0,X_k)$ are degenerated for any $k$.  This Markov chain can
also be represented (through an inversion of the time) as a dynamical
system $(T_{-n},\dots,T_{-1},T_0)$ which has the same law as
$(X_0,X_1,\dots,X_n)$ ($T_0$ and $X_0$ are random variables
distributed according to the invariant measure, see \cite{barb} for
more details). Let us recall the definition of a dynamical system.
\begin{definition}[dynamical system]
  A one-dimensional dynamical system is defined by
\begin{equation}\label{sysdyn}\forall
    k\in\mathbb{N}\,,\,T_k:=F^k(T_0)\;,
\end{equation}
where $F:I\to I$, $I$ is a compact subset of $\mathbb{R}$ and in this
context, $F^k$ denotes the $k$-th iterate of the appplication $F$:
$F^1=F$, $F^{k+1}= F \circ F^k$, $k\geq1$.  We assume that there
exists an invariant probability measure $\mu_0$, i.e.
$F(\mu_0)=\mu_0$, absolutely continuous with respect to Lebesgue's
measure, and that $T_0$ is a random variable with distribution
$\mu_0$.
\end{definition}

\noindent We  restrict our study  to one-dimensional dynamical
systems $T$ in the class $ {\cal F}$ of dynamical systems defined by a
transformation $F$ that satisfies the following assumptions (see
\cite{p}).
\begin{itemize}
\item $\forall k\in\mathbb{N}$, $\forall x\in$ int$(I)$, $\lim_{t\to
    0^+} F^k(x+t)=F^k(x^+)$ and $\lim_{t\to 0^-} F^k(x+t)=F^k(x^-)$
  exist;
\item $\forall k\in\mathbb{N}^*$, denoting $D^k_+=\{x\in$ int$(I),F^k(x^+)=x\}$
and $D^k_-=\{x\in$ int$(I),F^k(x^-)=x\}$, we assume
$\lambda\left(\displaystyle \bigcup_{k\in\mathbb{N}^*}\left(D^k_+\bigcup
D^k_-\right)\right)=0$, where $\lambda$ is the Lebesgue measure.

\end{itemize}

When the joint distributions of the pairs $(X_j,X_k)$ are not assumed
absolutely continuous (and then [H4] is not satisfied), we shall
instead assume:
\begin{itemize}
\item[] [H5] The dynamical system $(X_n)_{n\in\mathbb{Z}}$ belongs to
  $\mathcal F$.
\end{itemize}
\noindent We consider in this paper linear estimators as in
(\ref{estlin}).  The sequence of kernels $K_{m}$ is assumed to satisfy
the following assumptions.

\begin{description}
\item{(a)} The support of $K_{m}$ is a compact set with diameter
  ${O}(1/m^{1/d})$;
\item{(b)} The functions $x\mapsto K_{m}(x,y)$ and $x\mapsto
  K_{m}(y,x)$ are Lipschitz functions with Lipschitz constant
  ${O}\left(m^{1+1/d}\right)$;
\item{(c)} For all $x$ in the support of $K_{m}$, $\int K_{m}(x,y) \, \D y = 1$;
\item{(d)} The bias of the estimator $\hat f_n$ defined
  in~(\ref{estlin}) is of order $m_n^{-\rho/d}$, uniformly on compact
  sets.
\begin{gather}
  \sup_{\|x\| \leq M} \left| \mathbb{E} [\hat f_n(x)] - f(x) \right| =
  O (m_n^{-\rho/d} ) \; .
\end{gather}

\end{description}

\subsection{Results} \label{sec:results}
In all our results we consider kernels $K_m$ and a density estimator
of the form~(\ref{estlin}) such that assumptions (a), (b), (c) and
(d) hold.

\begin{theorem}[$\mathbb{L}^q$-convergence] \label{theorem:1}
\begin{description}
\item[\textbf{Geometric case.}] Under Assumptions [H4] or [H5] and
  [H1] or [H1'], the sequence $m_n$ can be chosen such that inequality
  (\ref{mini1}) holds for all $0<q<+\infty$.
\item[\textbf{Riemannian case.}] Under the assumptions [H4] or [H5],
  if additionally
\begin{itemize}
\item{}[H2] holds with $a > \max \left(1+2/d+(d+1)/\rho,2+1/d\right)$
  ($\eta$-dependence),
\item or [H2'] holds with $ a>1+2/d+1/\rho$ ($\tilde
  \phi$-dependence),
\end{itemize}
then the sequence $m_n$ can be chosen such that inequality
(\ref{mini1}) holds for all $0<q\le q_0=2\left\lceil
  (a-1)/2\right\rceil$.
\end{description}
\end{theorem}

\begin{theorem} [Uniform rates, geometric decays] \label{thorem:2}
  For any $M>0$, under Assumptions [H4] or [H5] and [H1] or [H1'] we
  have, for all $0<q<+\infty$, and for a suitable choice of the
  sequence $m_n$,
  \begin{eqnarray*}\mathbb{E}\sup_{\|x\|\le M}|\hat
    f_n(x)-f(x)|^q & = & O
       \left(\left(\frac{\log^{2(b+1)/b}(n)}{n}\right)^{q\rho/(d+2\rho)}\right)\;, \\
\sup_{\|x\|\le M}|\hat f_n(x)-f(x)|&=_{a.s.}&
O \left(\left(\frac{\log^{2(b+1)/b}(n)}{n}\right)^{\rho/(d+2\rho)}\right)\;.
\end{eqnarray*}
\end{theorem}

\begin{theorem}[Uniform rates, Riemannian decays] \label{theorem:3}
  For any $M>0$, under Assumptions [H4] or [H5], [H2] or [H2'] with
  $a\ge4$ and $\rho> 2d$, for $q_0 = 2 \lceil(a-1)/2 \rceil$ and $q\le
  q_0$, the sequence $m_n$ can be chosen such that
\begin{gather*}
  \mathbb{E}\sup_{\|x\|\le M} |\hat f_n(x)-f(x)|^q =
  O\left(n^{-\frac{q\rho }{d+2\rho+2d/(q_0+d)}}\right)\;,
\end{gather*}
or such that
\begin{gather*}
  \sup_{\|x\|\le M} | \hat f_n(x)-f(x)| =_{a.s.}
  O\left(\left(\frac{\log^{q_0+d}(n)}{n^{q_0-2} }\right)^{\frac{\rho}{d(q_0+2)+\rho(q_0+d)}}\right)\;.
\end{gather*}
\end{theorem}

{\bf Remarks.}
\begin{itemize}
\item Theorem~\ref{theorem:1} shows that the optimal convergence rate
  of (\ref{mini1}) still holds in the weak dependence context. In the
  Riemannian case, when $a\ge 4$, the conditions are satisfied if the
  density function $f$ is sufficient regular, namely, if $\rho> d+1$.

\item The loss with respect to the i.i.d. case in the uniform
  convergence rates (Theorems~\ref{thorem:2} and~\ref{theorem:3}) is
  due to the fact that the probability inequalities for dependent
  observations are not as good as Bernstein's inequality for i.i.d.
  random variables (Bernstein inequalities in weak dependence context
  are proved in \cite{neum}).  The convergence rates depend on the
  decay of the weak dependence coefficients. This is in contrast to
  the case of independent observations.

\item In Theorem~\ref{thorem:2} the loss is a power of the logarithm
  of the number of observations. Let us remark that this loss is
  reduced when $b$ tends to infinity. In the case of $\eta$-dependence
  and geometric decreasing, the same result is in \cite{dou2} for the
  special case $b=1$. In the framework of $\tilde \phi$-dependence,
  Theorem~\ref{thorem:2} seems to provide the first result on uniform
  rates of convergence for density estimators.

\item In Theorem~\ref{theorem:3}, the rate of convergence in the mean
  is better than the almost sure rate for technical reasons. Contrary
  to the geometric case, the loss is no longer logarithmic but is a
  power of $n$. The rate gets closer to the optimal rate as
  $q_0\to\infty$, or equivalently $a\to\infty$.

\item These results are new under the assumption of Riemannian decay
  of the weak dependence coefficients. The condition on $a$ is similar
  to the condition on $\beta$ in \cite{ang2}. Even if the rates are
  better than in \cite{dou3}, there is a huge loss with respect to the
  mixing case. It would be interesting to know the minimax rates of
  convergence in this framework.
\end{itemize}

\section{Models, applications and extensions}
\label{sec:applications}
The class of weak dependent processes is very large. We apply our
results to three examples: {\bf two-sided moving averages}, {\bf
  bilinear models} and {\bf expanding maps}. The first two will be
handled with the help of the coefficients $\eta$, the third one with
the coefficients $\tilde \phi$.

\subsection{Examples of $\eta$-dependent time series.}

\noindent
It is of course possible to define $\eta$-dependent random fields (see
\cite{ded3} for further details); for simplicity, we only consider
processes indexed by $\mathbb{Z}$.

\begin{definition}[Bernoulli shifts]
  Let $H:\mathbb{R}^{\mathbb{Z}}\to\mathbb{R}$ be a measurable
  function. A Bernoulli shift is defined as
  $X_n=H(\xi_{n-i},{i\in\mathbb{Z}})$ where $(\xi_i)_{i\in
    \mathbb{Z}}$ is a sequence of i.i.d random variables called the
  innovation process.
\end{definition}

\noindent
In order to obtain a bound for the coefficients $\{\eta_r\}$, we
introduce the following regularity condition on $H$. There exists a
sequence $\{\delta_r\}$ such that
$$
\sup_{i\in\mathbb{Z}}\mathbb{E}\left|H\left(\xi_{i-j},j\in\mathbb{Z}\right)-
  H\left(\xi_{i-j}\mathbbm{1}_{|j|<r},j\in\mathbb{Z}\right)\right|\le\delta_r
\;,
$$
{\bf Bernoulli shifts} are $\eta$-dependent with
$\eta_r=2\delta_{r/2}$ (see \cite{dou2}). In the following, we
consider two special cases of {\bf Bernoulli shifts}.

\begin{enumerate}
\item {\bf Non causal linear processes.} A real valued sequence
  $(a_i)_{i\in\mathbb{Z}}$ such that $\sum_{j\in\mathbb{Z}} a_j^2 <
  \infty$ and the innovation process $\{\xi_n\}$ define a {\bf
    non-causal linear process}
  $X_n=\sum_{-\infty}^{+\infty}a_i\xi_{n-i}$. If we control a moment
  of the innovations, the {\bf linear process} $(X_n)$ is
  $\eta$-dependent.  The sequence $\{\eta_r\}_{r\in\mathbb{N}}$ is
  directly linked to the coefficients $\{a_i\}_{i\in\mathbb{Z}}$ and
  various types of decay may occur.  We consider only Riemannian
  decays $\displaystyle a_i={\cal O}\left(i^{-A}\right)$ with $A\ge 5$
  since results for geometric decays are already known. Here $
  \eta_r={\cal O}\left(\sum_{|i|>r/2}a_i\right)={O}(r^{1-A})$ and [H2]
  holds. Furthermore, we assume that the sequence
  $(\xi_i)_{i\in\mathbb{Z}}$ is i.i.d. and satisfies the condition
  $|\mathbb{E} e^{iu\xi_0}|\le C(1+|u|)^{-\delta}$, for all
  $u\in\mathbb{R}$ and for some $\delta>0$ and $C<\infty$. Then, the
  densities $f$ and $f_{j,k}$ exist for all $j\ne k$ and they are
  uniformly bounded (see the proof in the causal case in Lemma 1 and
  Lemma 2 in \cite{gks}); hence [H4] holds. If the density $f$ of
  $X_0$ is $\rho$-regular with $\rho>2$, our estimators converge to
  the density with the rates:
\begin{itemize}
\item $n^{-\rho/(2\rho+1)}$ in $\mathbb{L}^q$-norm ($q\le 4$) at each
  point $x$,
\item $ n^{-\rho/(2\rho +3/2)}$ in $\mathbb{L}^q$-norm ($q\le 4$)
  uniformly on an interval,
\item $\left(\log^{4}(n)/n\right)^{\rho/(4\rho+3)}$
  almost surely on an interval.
\end{itemize}
In the first case, the rate we obtain is the same as in the i.i.d.
case. For such linear models, the density estimator also satisfies the
Central Limit Theorem (see \cite{hlt} and \cite{ded0}).

\item {\bf Bilinear model.} The process $\{X_t\}$ is a {\bf bilinear
    model} if there exist two sequences $(a_i)_{i\in\mathbb{N}^*}$ and
  $(b_i)_{i\in\mathbb{N}^*}$ of real numbers and  real
  numbers $a$ and $b$ such that:
\begin{equation} \label{bil}
X_t=\xi_t\left(a+\sum_{j=1}^\infty
a_jX_{t-j}\right)+b+\sum_{j=1}^\infty b_jX_{t-j}\; .
\end{equation}
Squared {\bf ARCH($\infty$)} or {\bf GARCH($p,q$)} processes satisfy
such an equation, with $b=b_j=0$ for all $j\geq1$.  Define
\begin{gather*}
  \lambda = \|\xi_0\|_{p} \sum_{j=1}^\infty a_j + \sum_{j=1}^\infty
  b_j \; .
\end{gather*}
If $\lambda < 1$, then the equation (\ref{bil}) has a strictly
stationary solution in $L^p$ (see \cite{dou4}). This solution is a
{\bf Bernoulli shift} for which we have the behavior of the
coefficient $\eta$:
\begin{itemize}
\item $\eta_r={O} \left(\E^{-\lambda r}\right)$ for some $\lambda>0$ if
  there exists an integer $N$ such that $a_i=b_i=0$ for $i\ge N$.
\item $\eta_r={O}(\E^{-\lambda\sqrt r})$ for some $\lambda>0$ if
  $a_i={O}(\E^{-Ai})$ and $b_i={O}(\E^{-B i})$ with $A>0$ and $B>0$.
\item $\eta_r={O} ( \left\{r/\log (r)\right\}^{-\lambda})$ for some
  $\lambda>0$ if $a_i={O}(i^{-A})$ and $b_i={O}(i^{-B})$ with $A>1$
  and $B>1$.
\end{itemize}

\noindent Let us assume that the i.i.d. sequence
$\{\xi_t\}$ has a marginal density $f_\xi\in C_\rho$, for some $\rho>
2$. The density of $X_t$ conditionally to the past can be written as a
function of $f_\xi$. We then check recursively that the common density
of $X_t$ for all $t$, say $f$, also belongs to $C_\rho$.  Furthermore,
the regularity of $f_\xi$ ensures that $f$ and the joint densities
$f_{j,k}$ for all $j\neq k$ are bounded (see \cite{dou4}) and [H4]
holds.  The assumptions of Theorem~\ref{theorem:1} are satisfied, and
the estimator $\hat f_n$ achieves the minimax bound (\ref{mini1}) if
either:
\begin{itemize}
\item There exists an integer $N$ such that $a_i=b_i=0$ for $i\ge N$;
\item There exist $A>0$ and $B>0$ such that $a_i={O}(e^{-Ai})$ and
  $b_i={O}(e^{-Bi})$;
\item There exist $A\ge 4$ and $B\ge 5$ such that $a_i={O}(i^{-A})$
  and $b_i={O}(i^{-B})$. Then, this optimal bound holds only for $2\le
  q< q(A,B)$ where $q(A,B)=2[((B-1)\wedge A)/2]$.
\end{itemize}
Note finally that the rates of uniform convergence provided by
Theorems~\ref{thorem:2} and~\ref{theorem:3} are sub-optimal.
\end{enumerate}

\subsection{Examples of $\tilde \phi$-dependent time series.}

Let us introduce an important class of {\bf dynamical
systems}:
\begin{example}
  $(T_i=F^i(T_0))_{i\in\mathbb{N}}$ is an {\bf expanding map} or
  equivalently $F$ is a Lasota-Yorke function if it satisfies the
  three following criteria.
\begin{itemize}
\item(Regularity) There exists a grid $0=a_0\le a_1\cdots\le a_n=1$
 such as $F\in \mathcal C_1$ and $|F'(x)|>0$ on $]a_{i-1},a_i[$ for each $i=1,\dots,n$.
\item(Expansivity) Let $I_n$ be the set on which $(F^n)'$ is defined.
  There exists $A>0$ and $s>1$ such that $\inf_{x\in I_n}|(F^n)'| > A
  s^n$.
\item(Topological mixing) For any nonempty open sets $U$, $V$, there
exists $n_0\ge 1$ such as $F^{-n}(U)\cap V\neq\varnothing$ for all
$n\ge n_0$.
\end{itemize}
\end{example}

\noindent Examples of {\bf Markov chains}
$X_n=G(X_{n+1},\epsilon_n)$ associated to an {\bf expanding map}
$\{T_n\}$ belonging to ${\cal F}$ are given in \cite{barb} and
\cite{ded2}. The simplest one is
$X_k=\left(X_{k-1}+\epsilon_k\right)/2$ where the $\epsilon_k$ follows
a binomial law and $X_0$ is uniformly distributed on $[0,1]$. We
easily check that $F(x)=2x$ mod$\,1$, the transformation of the
associated {\bf dynamical system} $T_n$, satisfies all the assumptions
such as $T_n$ is an {\bf expanding map}
belonging to ${\cal F}$.\\

\noindent The  coefficients of $\tilde\phi$-dependence of such a  {\bf Markov
  chain} satisfy $\tilde \phi(r)=O(e^{-ar})$ for some $a>0$ (see
\cite{ded2}).  Theorems~\ref{theorem:1} and~\ref{thorem:2} give the
$\mathbb{L}^q$ rate $ n^{-\rho/(2\rho+1)}$, the uniform $\mathbb{L}^q$
rate and the almost sure rate $ \left( \log^4 (n)/n
\right)^{\rho/(2\rho+1)}$ of the estimators of the density of $\mu_0$.

\subsection{Sampled process}

Since we do not assume stationarity of the observed process, the
following observation scheme is covered by our results.  Let
$(x_n)_{n\in\mathbb{Z}}$ be a stationary process whose marginal
distribution is absolutely continuous, let $(h_n)_{n\in\mathbb{Z}}$ be
a sequence of monotone functions and consider the sampled process
$\{X_{i,n}\}_{1 \le i \le n}$ defined by $X_{i,n} = x_{h_n(i)}$. The
dependence coefficients of the sampled process may decay to zero
faster than the underlying unoberved process. For instance, if the
dependence coefficients of the process $(x_n)_{n\in\mathbb{Z}}$ have a
Riemannian decay, those of the sampled process $\{x_{h_n(i)}\}$ with
$h_n(i)=i2^n$ decay geometrically fast.  The observation scheme is
thus a crucial factor that determines the rate of convergence of
density estimators.

\subsection{Density estimators and bias}
\label{sec:bias}

In this section, we provide examples of kernels $K_{m}$ and
smoothness assumptions on the density $f$ such that assumptions (a),
(b), (c) and (d) of subsection~\ref{nes} are satisfied.

\noindent {\bf Kernel estimators} The kernel estimator associated to
the bandwidth parameter $m_n$ is defined by:
$$
\hat{f}_{n}(x)=\frac{m_n}{n}\sum_{i=1}^nK\left(m_n^{1/d}(x-X_i)\right)\;.
$$
We briefly recall the classical analysis for the deterministic part
$R_n$ in this case (see \cite{t}). Since the sequence $\{X_n\}$ has
a constant marginal distribution, we have $\mathbb{E}[\hat f_n(x)] =
f_n(x)$ with $f_n(x)=\int_DK(s)f\left(x-s/m_n^{1/d}\right) \D s$.
Let us assume that $K$ is a Lipschitz function compactly supported
in $D\subset\mathbb{R}^d$.  For $\rho>0$, let $K$ satisfy, for all
$j=j_1+\cdots+j_d$ with $(j_1,\dots,j_d)\in\mathbb{N}^d$:
\begin{equation*}
\int x_1^{j_1} \cdots x_d^{j_d} K(x_1,\ldots,x_d) \D x_1\cdots \D x_d=
\begin{cases} 1 & \text{if}\,j=0,\\ 0 & \text{for }\,j\in\{1,\dots,\lceil \rho-1\rceil-1\},\\
\neq 0 & \ \text{if}\ j=\lceil \rho-1\rceil.
\end{cases}\end{equation*}
Then the kernels $K_{m}(x,y)=mK\left(m^{1/d}(x-y)\right)$ satisfy
(a), (b) and (c). Assumption (d) holds and if $f\in\mathcal C_\rho$,
where $C_\rho$ is the class of function $f$ such that for
$\rho=\lceil \rho-1 \rceil+c$ with $0<c\le 1$, $f$ is $\lceil \rho-1
\rceil$-times continuously differentiable and there exists $A>0$
such that $\forall(x,y) \in \mathbb{R}^d\times\mathbb{R}^d$,
$|f^{(\lceil \rho-1
  \rceil)}(x)-f^{(\lceil \rho-1 \rceil)}(y)|\le A|x-y|^c$.

\medskip

\noindent {\bf Projection estimators} We only consider in this
section the case $d=1$. Under the assumption that the family
$\{1,x,x^2,\dots\}$ belongs to $L^2(I,\mu)$, where $I$ is a bounded
interval of $\mathbb{R}$ and $\mu$ is a measure on $I$, an
orthonormal basis of $L^2(I,\mu)$ can be defined which consists of
polynomials $\{P_0, P_1,P_2,\dots\}$.   We assume that $f$ belongs
to a class $\mathcal C'_\rho$ which is slightly more restrictive
than the class $\mathcal C_\rho$ (see Theorem 6.23 p.218 in
\cite{dou5} for details). Then for any $f\in L^2(I,\mu)\cap\mathcal
C'_\rho$, there exists a function $\pi_{f,m_n}\in V_{m_n}$ such that
$\sup_{x\in
  I}|f(x)-\pi_{f,m_n}(x)|={O}(m_n^{-\rho})$.  Consider then the
projection $\pi_{m_n}f$ of $f$ on the subspace
$V_{m_n}=$Vect$\{P_0,P_1,\dots,P_{m_n}\}$. It can be expressed as
$$
\pi_{m_n}f(x) = \sum_{j=0}^{m_n} \left\{ \int_I P_j(s)f(s) \D \mu(s)
\right\} P_j(x).
$$
The projection estimator of the density $f$ of the real valued
random variables $\{X_i\}_{1\le i\le n}$ is naturally defined as
\begin{equation*}\label{proj}
\hat f_n(x) = \frac{1}{n} \sum_{i=1}^n K_{m_n}(x,X_i)
= \frac{1}{n} \sum_{i=1}^n \sum_{j=0}^{m_n} P_j(X_i) P_j(x)\;.
\end{equation*}
Then $\displaystyle\mathbb{E}\hat f_n(x)=\pi_{m_n}f(x)$ is an
approximation of $f(x)$ in $V_{m_n}$. The fact that $I$ is compact
and the Christoffel-Darboux formula and its corollary (see
\cite{sz}) ensure properties (a) and (b) for the kernels $K_{m}$. We
easily check that properties (c) also holds. Unfortunately, the
optimal rate $(m_n^{-\rho})$ does not necessarily hold. We then have
to consider the weighted kernels $K^a_m(x,y)$ defined by:
$$
K^a_m(x,y)=\sum_{j=0}^{m}a_{m,j}\sum_{k=0}^jP_k(x)P_k(y)\;,
$$
where $\{a_{m,j};\,m\in\mathbb{N},\,0\le j\le m\}$ is a weight
sequence satisfying $ \sum_{j=0}^m a_{m,j}=1$ and for all $j$:
$\lim_{m\to\infty} a_{m,j}= 0$. If the sequence $\{a_{m,j}\}$ is
such that $K^a_m$ is a nonnegative kernel then $\|K^a_m\|_1=\int_I
K^a_m(x,s)d\mu(s)=1$ and the kernel $K^a_{m}$ satisfies (a), (b) and
(c). Moreover, the uniform norm of the operator $f \mapsto
K^a_m*f(x)$ is $\sup_{\|f\|_\infty=1} \|K^a_m*f\|_\infty =
\|K^a_m\|_1=1$.  The linear estimator built with this kernel is
$$
\hat f^a_n(x) = \frac{1}{n} \sum_{i=1}^n \sum_{j=0}^{m_n}
a_{{m_n},j} \sum_{k=0}^j P_k(X_i) P_k(x) \; ,
$$
and its bias has the optimal rate:
\begin{eqnarray*}
|\mathbb{E}\hat
f^a_n(x)-f(x)|&= &|K^a_{m_n}*f(x)-\pi_{f,m_n}f(x)+\pi_{f,m_n}f(x)-f(x)|\;,\\
&\le& |K^a_{m_n}*(f(x)-\pi_{f,m_n}f(x))+\pi_{f,m_n}f(x)-f(x)|\;,\\
&\le& (\|K^a_{m_n}\|_1+1)m_n^{-\rho}={\cal
O}(m_n^{-\rho})\;.
\end{eqnarray*}
Such an array $\{a_{m,j}\}$ cannot always be defined. We give an
example where it is possible.

\begin{example}[Fejer kernel]
  For the trigonometric basis $\{\cos
  (nx),\sin(nx)\}_{n\in\mathbb{N}}$, we can find a $2\pi$-periodic
  function $f\in \mathcal C'_1$ such that
  $\sup_{x\in[-\pi;\pi]}|f(x)-\pi_mf(x)|={O}(m^{-1}\log m)$. The
  associated estimator reads:
  $$
  \hat f_n(x) = \frac1{2\pi} + \frac1{n\pi} \sum_{i=1}^n
  \sum_{k=1}^{m_n} \cos(kX_i) \cos(kx) + \sin(kX_i) \sin(kx) \;.
  $$
  We remark that $\mathbb{E} \hat f_n$ is the Fourier series of $f$
  truncated at order $m_n$:
  $$
  D_{m_n}f(x) = \frac1{2\pi} \int_0^{2\pi} f(t) D_{m_n}(x-t) \D t
  \; .
$$
where
$$
D_{m}(x) = \sum_{k=-{m}}^{m} \E^{\I kx} =
\frac{\sin(\{2{m}+1\}x/2)}{\sin(x/2)}
$$
is (the symmetric) Dirichlet's kernel.  Recall that Fejer's kernel
is defined as
$$
F_{m}(x) = \frac1{m} \sum_{k=0}^{m-1} D_k(x) =
\sum_{k=-(m-1)}^{m-1} \left(1-\frac{|k|}{m}\right) \E^{\I kx} =
\frac{\sin^2({m}x/{2})}{m \, \sin^2(x/2)} \;.
$$
The kernel $F_m$ is a nonnegative weighted kernel corresponding to
Dirichlet's kernel and the sequence of weights $a_{m,j}=1/m$ and
satisfies (a), (b) and (c).  The estimator associated to the Fejer's
kernels is defined by
$$
\tilde f_n(x) = \frac1{2\pi} + \frac1{n\pi} \sum_{i=1}^n
\sum_{j=1}^{m_n} \frac{1}{m_n} \sum_{k=1}^j\cos k X_i \cos kx + \sin
kX_i\sin kx\;,
$$
If the common density $f$ is $2\pi$-periodic and belongs to
$\mathcal C'_1$, then assumption (d) holds.
\end{example}
Using general Jackson's kernels (see \cite{dou5}), we can find an
estimator such that $R_n={O}(m_n^{-\rho/d})$ for other values of
$\rho$, but the weight sequence $a_{m,j}$ highly depends of the value
of $\rho$.

\noindent {\bf Wavelet estimation} Wavelet estimation is a
particular case of projection estimation. For the sake of simplicity,
we restrict hte study to $d=1$.

\begin{definition} [Scaling function  \cite{d}]
  A function $\phi\in L^2(\mathbb{R})$ is called a scaling function if
  the family $\{\phi(\cdot - k)\,;\,k\in\mathbb{Z}\}$ is orthonormal.
\end{definition}

\noindent We  choose the bandwidth parameter
$m_n=2^{j(n)}$ and define $V_j = \mathrm{Vect} \{\phi_{j,k}, k \in
\mathbb{Z}\}$, where $\phi_{j,k}=2^{j/2}\phi(2^j(x-k))$. Under the
assumption that $\phi$ is compactly supported, we define (the sum over
the index $k$ is in fact finite):
$$
\hat{f}_{n}(x) = \frac{1}{n} \sum_{k=-\infty}^\infty \sum_{i=1}^n
\phi_{j(n),k}(X_i) \phi_{j(n),k}(x)\;.
$$
The wavelets estimator is of the form (\ref{estlin}) with $
K(x,y)=\sum_{k=-\infty}^\infty\phi(y-k)\phi(x-k)$ and
$K_m(x,y)=mK(mx,my)$. Under the additionnal assumption that
$\sum_{k\in\mathbb{Z}}\phi(x-k)= 1$ for almost all $x$, we can write:
\begin{eqnarray*}
\left|\mathbb{E}(\hat{f}_n(x)-f(x))\right|&\le&\left|\int K_{m_n}(y,x)f(y)dy-f(x)\right|\;,\\
&=&\left|\int {m_n} K({m_n}y,{m_n} x)(f(y)-f(x))dy\right|\;,\\
&=&\left|\int {m_n}K({m_n}x+t,{m_n}
x)(f(x+t/m_n)-f(x))dt\right|\;.
\end{eqnarray*}
If $\phi$ is a Lipschitz function such that $\int \phi(x)x^jdx=0$ if
$0<j<\lceil \rho-1\rceil$ and $\int \phi(x)x^{\lceil
\rho-1\rceil}dx\neq 0$, then the kernel $K_m$ satisfy properties
(a), (b) and (c). If $f\in C_\rho$, then Assumption (d) holds.

\section{Proof of the Theorems} \label{sec:proof}

The proof of our results is based on the decomposition:
\begin{equation}\label{fb}
\hat{f}_n(x)-f(x)=\underbrace{\hat{
f}_n(x)-\mathbb{E}\left(\hat{f}_n(x)\right)}_{F\!L_{n}(x) =
\text{fluctuation}}
+\underbrace{\mathbb{E}\left(\hat{f}_n(x)\right)-f(x)}_{
\text{bias}}\;.
\end{equation}
The bias term is of order $m_n^{-\rho/d}$ by Assumption (d).
We now present three lemmas useful to derive the rate of the
fluctuation term.

\begin{lemma}[Moment inequalities] \label{imp}
  For each even integer $q$, under the assumption [H4] or [H5] and if
  moreover one of the following assumption holds:
\begin{itemize}
\item{} [H1] or [H1'] holds (geometric case);
\item{} [H2] holds, $m_n=n^{\delta}\log(n)^\gamma$ with $\delta>0$,
  $\gamma \in \mathbb{R}$ and
  $$
  a>\max\left(q-1,\frac{(q-1)\delta(4+2/d)}{q-2+\delta(4-q)},2+\frac{1}{d}\right)\;,
  $$
\item{} [H2'] holds, $m_n=n^{\delta}\log(n)^\gamma$ with $\delta>0$ and
  $\gamma \in \mathbb{R}$ and
  $$
  a > \max \left(q-1,\frac{(q-1) \delta(2+2/d)}{q-2+\delta(4-q)}, 1
    + \frac{1}{d}\right)\;.
  $$
\end{itemize}

Then, for each $x\in\mathbb{R}^d$,
$$
\limsup_{n\to\infty}  \; (n/m_n)^{q/2} \, \| F\!L_n(x)\|^q_q < +
\infty\;.
$$
\end{lemma}

\begin{lemma}[Probability inequalities]\label{liii}\

\begin{itemize}
\item Geometric case. Under Assumptions [H4] or [H5] and [H1] or [H1']
  there exist positive constants $C_1,C_2$ such that
  $$
  \mathbb{P}\left(|F\!L_n(x)|\ge \epsilon\sqrt{m_n/n} \right) C_1
  \leq \exp \{ - C_2\epsilon^{b/(b+1)} \} \;.
$$

\item Riemannian case. Under Assumptions [H4] or [H5], if $m_n =
  n^\delta \log(n)^\gamma$ and if one of the following assumtions
  holds:
\begin{itemize}
\item{} [H2] with $ a > \max \{1+2(\delta+1/d)/(1-\delta),2+1/d\}$,
\item{} [H2'] with $a > \max \left(1+2\{1/d(1-\delta)\}, 1+1/d
  \right)$,
\end{itemize}
then,
$$
\mathbb{P}\left(|F\!L_n(x)| \ge \epsilon \sqrt{m_n/n}\right) \leq C
\epsilon^{-q_0} \; ,
$$
with  $q_0=2\left\lceil (a-1)/2\right\rceil$.
\end{itemize}
\end{lemma}

\begin{lemma}[Fluctuation rates] \label{rat}
  Under the assumptions of Lemma \ref{liii}, we have for any $M>0$,
\begin{itemize}
\item Geometric case.
  $$
  \sup_{\|x\|\le M}|F\!L_n(x)| =_{a.s.} O\left(
    \sqrt{\frac{m_n}{n}}\log^{(b+1)/b}(n) \right) \; ;
  $$
\item Riemannian case.
$$\displaystyle\sup_{\|x\|\le M}|F\!L_n(x)|=_{a.s.}{ O}\left(
\left(\frac{m_n^{1+2/q_0}}{n^{1-2/q_0}}\right)^{\frac1{1+d/q_0}}\log n\right)\;,
$$
\end{itemize}
with $q_0=2\left\lceil (a-1)/2 \right\rceil$.
\end{lemma}

\noindent{\bf Remarks.}
\begin{itemize}
\item In Lemma \ref{imp}, we improve the moment inequality of
  \cite{dou3}, where the condition in the case of coefficient $\eta$
  is $a>3(q-1)$, which is always stronger than our condition.
\item In the i.i.d. case a Bernstein type inequality is available:
  $$
  \mathbb{P}\left(|F\!L_n(x)| \geq \epsilon \sqrt{\frac{m}{n}}
  \right) \leq C_1 \exp\left( -C_2 \epsilon^2 \right) \;,
  $$
  Lemma \ref{liii} provides a weaker inequality for dependent
  sequences. Other probability inequalities for dependent sequences
  are presented in \cite{ded2} and \cite{neum}.
  \item Lemma \ref{rat} gives the almost sure bounds for the
  fluctuation. It is derived directly from the two previous lemmas.
\end{itemize}

\subsubsection{Proof of the lemmas}
\begin{proof}[Proof of Lemma~\ref{imp}]
  Let $x$ be a fixed point in $\mathbb{R}^d$. Denote
  $Z_i=u_n(X_i)-\mathbb{E} u_n(X_i)$ where $u_n(.)=
  K_{m_n}(.,x)/\sqrt{m_n}$.  Then
\begin{equation}\label{id}
\sum_{i=1}^nZ_i=\sum_{i=1}^nu_n(X_i)-\mathbb{E}
u_n(X_i)=\frac{n}{\sqrt{m_n}}(\hat f_n(x)-\mathbb{E} \hat
f_n(x))=\frac{n}{\sqrt{m_n}}F\!L_{n}(x)\;.
\end{equation}
The order of magnitude of the fluctuation $FL_n(x)$ is obtained by
applying the inequality (\ref{douklou}) to the centered sequence
$\{Z_i\}_{1\le i\le n}$ defined above. We then control the normalized
fluctuation of (\ref{id}) with the covariance terms $C_k(r)$ defined
in equation (\ref{ckr}). Firstly, we bound the covariance terms:
\begin{itemize}
\item{\bf Case $r=0$.} Here $t_1=\dots=t_k=i$. Then we get:
$$C_k(r)=\left|\mathrm{cov}\left(Z_{t_1}\cdots Z_{t_p},
 Z_{t_{p+1}}\cdots Z_{t_k}\right)\right|\le 2\mathbb{E}|Z_i|^k\;.$$ By
 definition of $Z_i$:
\begin{equation}\label{c0}
\mathbb{E}|Z_i|^k\le 2^k \mathbb{E}|u_n(X_i)|^k\le 2^k
\|u_n\|_\infty^{k-1}\mathbb{E}|u_n(X_0)|\;.
\end{equation}
\item{\bf Case  $r>0$.} $\displaystyle
C_k(r)=\left|\mathrm{cov}\left(Z_{t_1}\cdots Z_{t_p},
 Z_{t_{p+1}}\cdots Z_{t_k}\right)\right|$ is
bounded in different ways, either using weak-dependence property or
by direct bound.
\begin{itemize}
\item{\bf Weak-dependence bounds:}
\begin{itemize}
\item{\it $\eta$-dependence:} Consider the following
application:
\begin{eqnarray*}
\phi_p:(x_1,\dots,x_p)\mapsto (u_n(x_1)\cdots u_n(x_p))\;.
\end{eqnarray*}
Then $\|\phi_p\|_\infty\le2^p\|u_n\|_\infty^p$ and $\mathrm{Lip}\, \phi_p\le
2^p\|u_n\|_\infty^{p-1}\mathrm{Lip}\, u_n$. Thus by $\eta$-dependence, for all
$k\ge 2$ we have:
\begin{eqnarray} \nonumber
C_k(r)&\le&\left(p2^p\|u_n\|_\infty^{p-1}+\right.\left.(k-p)2^{p-k}\|u_n\|_\infty^{p-k-1}\right)\mathrm{Lip}\,
u_n\eta_r\;,\\
\label{cre} & \le &k2^k\|u_n\|_\infty^{k-1}\mathrm{Lip}\, u_n\eta_r\;.
\end{eqnarray}

\item{\it $\tilde \phi$-dependence:} We use the inequality
  (\ref{tau}). Using the bound
$$\mathbb{E}|\phi_p(X_1,\dots,X_p)|\le
\|u_n\|_\infty^{p-1}\mathbb{E}|u_n(X_0)|\;,$$ we derive a bound for the
covariance terms:
\begin{eqnarray}\label{ctau}
C_k(r)&\le& k2^k\|u_n\|_\infty^{k-2}\mathbb{E}|u_n(X_0)|\mathrm{Lip}\, u_n\tilde
\phi(r)\;.
\end{eqnarray}
\end{itemize}

\item{\bf Direct bound:} Triangular inequality implies for $C_k(r)$:
$$\left|\mathrm{cov}\left(Z_{t_1}\cdots Z_{t_p},Z_{t_{p+1}}\cdots Z_{t_k}\right)\right|
\le\underbrace{\left|\mathbb{E}\prod_{i=1}^k
Z_{t_i}\right|}_A+\underbrace{\left|\mathbb{E}\prod_{i=1}^p
Z_{t_i}\right|}_{B_p}\underbrace{\left|\mathbb{E}\prod_{i=p+1}^k
Z_{t_i}\right|}_{B_{k-p}}\;,$$
\begin{eqnarray*}
A &=&\left|\mathbb{E}\left(u_n(X_{t_1})-\mathbb{E}
u_n(X_{t_1})\right)\cdots\left(u_n(X_{t_k})-\mathbb{E} u_n(X_{t_k})
\right)\right|\;,\\
&=&\left|\mathbb{E} u_n(X_0)\right|^{k}+\left|\mathbb{E}\left(u_n(X_{t_1})\cdots
u_n(X_{t_k})\right)\right|
\\&&+\sum_{s=1}^{k-1}\left|\mathbb{E} u_n(X_0)\right|^{k-s}\sum_{t_{i_1}\le\cdots\le
t_{i_s}}\left|\mathbb{E}\left(u_n(X_{t_{i_1}})\cdots
u_n(X_{t_{i_s}})\right)\right|\;.
\end{eqnarray*}
Firstly, with $k\ge 2$:
\begin{eqnarray*} \left|\mathbb{E} u_n(X_0)\right|^{k}\le
\|u_n\|_\infty^{k-2}(\mathbb{E}|u_n(X_0)|)^2\;.\end{eqnarray*} Secondly, if
$1\le s\le k-1$:
\begin{eqnarray*}
\left|\mathbb{E}\left(u_n(X_{t_{i_1}})\cdots
u_n(X_{t_{i_s}})\right)\right|&\le&\mathbb{E}|u_n(X_{t_{i_1}})\cdots
u_n(X_{t_{i_s}})|\;,\\&\le&\|u_n\|_{\infty}^{s-1}\mathbb{E}|u_n(X_0)|\;,\\
\left|\mathbb{E} u_n(X_0)\right|^{k-s}&\le&
\|u_n\|_\infty^{k-s-1}\mathbb{E}|u_n(X_0)|\;.\end{eqnarray*} Thirdly there
is at least two different observations with a gap of $r>0$ among
$X_{t_1},\dots,X_{t_k}$ so for any integer $k \ge 2$
:\begin{eqnarray*} \left|\mathbb{E}\left(u_n(X_{t_1})\cdots
u_n(X_{t_k})\right)\right|\le\|u_n\|_{\infty}^{k-2}\mathbb{E}|u_n(X_0)u_n(X_r)|\;.
\end{eqnarray*}
Then, collecting the last four inequations yields:
\begin{eqnarray*}  A&\le&\|u_n\|_\infty^{k-2}(\mathbb{E}|u_n(X_0)|)^{2}
\\
&&+(\mathbb{E}|u_n(X_0)|)^2\sum_{s=1}^{k-1}C_s^{k}\|u_n(X_0)\|_\infty^{k-2}
+\|u_n\|_{\infty}^{k-2}\mathbb{E}|u_n(X_0)u_n(X_r)|\;.
\end{eqnarray*}
So:\begin{equation}\label{AA} A\le
\|u_n\|_\infty^{k-2}\left((2^k-1)(\mathbb{E}|u_n(X_0)|)^{2}+
\mathbb{E}|u_n(X_0)u_n(X_r)|\right)\;.
\end{equation}
Now, we bound $B_i$ with $i<k$. As before:
\begin{eqnarray*}
B_i&=&\left|\mathbb{E}\left(u_n(X_{t_1})-\mathbb{E} u_n(X_{t_1})\right)\cdots\left(u_n(X_{t_i})-\mathbb{E} u_n(X_{t_i})\right)\right|\;,\\
&=&\sum_{s=0}^{i}\left|\mathbb{E}(u_n(X_0)\right|^{i-s}\sum_{t_{j_1}\le\dots\le
t_{j_s}}\left|\mathbb{E}\left(u_n(X_{t_{j_1}})\cdots
u_n(X_{t_{j_s}})\right)\right|\;,\\
&\le& 2^i\|u_n\|_\infty^{i-2}(\mathbb{E}|u_n(X_0)|)^2\;.
\end{eqnarray*}
Then: \begin{equation}\label{BB} B_p\times B_{k-p}\le
2^k\|u_n\|_\infty^{k-4}(\mathbb{E}|u_n(X_0)|)^4\le2^k\|u_n\|_\infty^{k-2}(\mathbb{E}|u_n(X_0)|)^2\;.
\end{equation}
Another interesting bound for $r>0$ follows, because according to
inequalities (\ref{AA}) and (\ref{BB}) we have:
$$C_{k}(r) \le \|u_n\|_\infty^{k-2}\left((2^{k+1}-1)(\mathbb{E}|u_n(X_0)|)^{2}+
\mathbb{E}|u_n(X_0)u_n(X_r)| \right)\;.$$ Noting
$\gamma_n(r)=\mathbb{E}|u_n(X_0)u_n(X_r)|\vee (\mathbb{E}|u_n(X_0)|)^2$, we have:
\begin{equation}\label{crd}
C_k(r)\le 2^{k+1}\|u_n\|_\infty^{k-2}\gamma_n(r)\;.
\end{equation}
\end{itemize}
\end{itemize}
We now use the different values of the bounds in inequalities
(\ref{c0}), (\ref{cre}), (\ref{ctau}) and (\ref{crd}). If we define
the sequence $(w_r)_{0\le r\le n-1}$ as:
\begin{itemize}
\item $w_0=1$, \item $\displaystyle
w_r=\gamma_n(r)\wedge\|u_n\|_\infty \mathrm{Lip}\, u_n
\eta_r\wedge\mathbb{E}|u_n(X_0)|\mathrm{Lip}\, u_n\tilde \phi(r)$,
\end{itemize}
then, for all $r$ such that $0\le r\le n-1$ and for all $k\ge 2$:
$$
\displaystyle C_k(r)\le k2^k\|u_n\|_\infty^{k-2}w_r\;.
$$
We derive from this inequality and from (\ref{douklou}):
\begin{eqnarray*}
\left\|\sum_{i=1}^nZ_i\right\|_q^q&\le&
\frac{(2q-2)!}{(q-1)!}\left\{\left(n\sum_{r=0}^{n-1}C_2(r)\right)^{q/2}\vee
n\sum_{r=0}^{n-1}(r+1)^{q-2}C_q(r)\right\}\;,\\
&\preceq&\left(q\sqrt{n}\right)^q\left\{\left(\sum_{r=0}^{n-1}w_r\right)^{q/2}\vee
\left(\frac{\|u_n\|_\infty}{\sqrt n}\right)^{q-2}
\sum_{r=0}^{n-1}(r+1)^{q-2}w_r\right\}\;.
\end{eqnarray*}
The symbol $\preceq$ means $\le$ up to an universal constant. In
order to control $w_r$, we give bounds for the terms
$\gamma_n(r)=\mathbb{E}|u_n(X_0)u_n(X_r)|\vee (\mathbb{E}|u_n(X_0)|)^2$:
\begin{itemize}
\item In the case of [H4], we have:
\begin{eqnarray*}&\mathbb{E}|u_n(X_0)u_n(X_r)|\le
\sup_{j,k}\|f_{j,k}\|_\infty\|u_n\|_1^2\;,\\&(\mathbb{E}|u_n(X_0)|)^2\le
\|f\|_\infty^2\|u_n\|_1^2\;.\end{eqnarray*}
\item In the case of [H5], Lemma 2.3 of \cite{p} proves that
  $\mathbb{E}|u_n(X_0)u_n(X_r)|\le (\mathbb{E}|u_n(X_0)|)^2$ for $n$
  sufficiently large and the same bound as above remains true for the
  last term.
\end{itemize}
In both cases, we conclude that $\gamma_n(r)\preceq \|u_n\|_1^2$.
The properties (a), (b) and (c) of section \ref{nes} ensures that
$\displaystyle\|u_n\|_1^2\preceq \frac{1}{m_n}$,
$\displaystyle\|u_n\|_\infty\mathrm{Lip}\, u_n\preceq m_n^{1+1/d}$
and $\displaystyle\mathbb{E}|u_n(X_0)|\mathrm{Lip}\, u_n\preceq
m_n^{1/d}$. We then have for $r\ge 1$:
\begin{equation}\label{beta}
w_r\preceq \frac{1}{m_n}\wedge m_n^{1+1/d}\eta_r\wedge
m_n^{1/d}\tilde \phi_r\;.
\end{equation}
In order to prove Lemma \ref{imp}, it remains to control the sums
\begin{equation}\label{sum}\left(\frac{\|u_n\|_\infty}{\sqrt
      n}\right)^{k-2}\sum_{r=0}^{n-1}(r+1)^{k-2}w_r\;,\end{equation}
for $k=2$ and $k=q$ in both  Riemannian and geometric cases.
\begin{itemize}
\item{\bf Geometric case.\\}
{\it Under [H1] or [H1']:} We remark that $a\wedge b\le a^\alpha
b^{1-\alpha}$ for all $\alpha \in [0;1]$. Using (\ref{beta}), we
obtain first that $w_r\preceq (\eta_r\wedge\tilde \phi_r)^\alpha
m_n^{\alpha(1+1/d)-(1-\alpha)}$ for $n$ sufficiently large. Then for
$0<\alpha\le\frac{d}{2d+1}$ we bound $w_r$ independently of $m_n$:
$w_r\preceq (\eta_r\wedge\tilde \phi_r)^\alpha$. For all even
integer $k\ge2$ we derive from the form of $\eta_r\wedge\tilde
\phi_r$ that (in the third inequality $u={a}r^b$):
\begin{eqnarray*}
\sum_{r=1}^{n-1}(r+1)^{k-2}w_r&\preceq&\sum_{r=0}^{n-1}(r+1)^{k-2}\exp(-\alpha ar^b)\;,\\
&\preceq&\int_0^\infty r^{k-2}\exp(-\alpha ar^b)dr\;,\\
&\preceq&\frac1{ba^{\frac{k-1}b} }\int_1^\infty
u^{\frac{k-1}b-1}\exp(-u)du\;,\\
&\preceq&\frac1{ba^{\frac{k-1}b} }
\Gamma\left(\frac{k-1}{b}\right)\;.
\end{eqnarray*}
Using the Stirling formula, we can find a constant $B$ such that, for
the special cases $k=2$ and $k=q$:
$$
\sum_{r=1}^{n-1}(r+1)^{k-2}w_r\preceq\frac1{ba^{\frac{k-1}b} }
\Gamma\left(\frac{k-1}{b}\right)\preceq (Bk)^{\frac kb}\;.
$$

\item{\bf Riemannian case.\\} {\it Under [H6] and [H2]:} Let us
  recall that [H6] implies that $m_n\le n^\delta$ for $n$ sufficiently
  large and $0<\delta<1$ and that the assumption of Lemma \ref{imp}
  implies that:
  $$\displaystyle
  a>\max\left(q-1,\frac{\delta(q-1)(4+2/d)}{q-2+\delta(4-q)},2+\frac{1}{d}\right)\;.$$
  Then, we have $\displaystyle
  a>\max\left(k-1,\frac{\delta(k-1)(4+2/d)}{k-2+\delta(4-k)}\right)$
  for both cases $k=q$ or $k=2$. This assumption on $a$ implies that:
  $$
  \frac{(k+2/d)\delta+2-k}{2(a-k+1)}<\frac{(4-k)\delta+k-2}{2(k-1)}\;.
  $$
  Furthermore, reminding that $0<\delta<1$:
  $$
  0<\frac{(4-k)\delta+k-2}{2(k-1)}=1-\frac{k(1+\delta)-4\delta}{2(k-1)}\le1\;.
  $$
  We derive from the two previous inequalities that there exists
  $\zeta_k\in]0,1[$ verifying $\displaystyle
  \frac{(k+2/d)\delta+2-k}{2(a-k+1)} < \zeta_k <
  \frac{(4-k)\delta+k-2}{2(k-1)}.$

\noindent For $k=q$ or $k=2$, we now use Tran's technique as in
\cite{ang1}. We divide the sum (\ref{sum}) in two parts in order to
bound it by sequences tending to $0$, due to the choice of $\zeta_k$:
\begin{eqnarray*}
\left(\sqrt\frac{{m_n}}{n}\right)^{k-2}\sum_{r=0}^{[n^{\zeta_k}]-1}(r+1)^{k-2}w_r&\preceq&
\left(\sqrt\frac{{m_n}}{n}\right)^{k-2}\frac{[n^{\zeta_k}]^{k-1}}{{m_n}}\;,\\
&\preceq&n^{(2\zeta_k(k-1)-((4-k)\delta+k-2))/2}\;,\\&=&{O}(1)\;,\\
\left(\sqrt\frac{{m_n}}{n}\right)^{k-2}\sum_{r=[n^{\zeta_k}]}^{n-1}(r+1)^{k-2}w_r&\le&
\left(\sqrt\frac{{m_n}}{n}\right)^{k-2}m_n^{1+1/d}[n^{\zeta_k}]^{k-1-a}\;,\\
&\le&n^{(-2\zeta_k(a-k-1)+((k+2/d)\delta+2-k))/2}\;,\\
 &=&{O}(1)\;.
\end{eqnarray*}
{\it Under [H6] and [H2']:} Under the assumption of Lemma \ref{imp}:
$$\displaystyle
a>\max\left(q-1,\frac{\delta(q-1)(2+2/d)}{q-2+\delta(4-q)},1+\frac{1}{d}\right)\;,$$
we derive exactly as in the previous case that there exists
$\zeta_k\in]0;1[$ for $k=q$ or $k=2$ such that
$$\displaystyle
\frac{(k-2+2/d)\delta+2-k}{2(a-k+1)}<\zeta_k<\frac{(4-k)\delta+k-2}{2(k-1)}\;.$$
We then apply again the Tran's technique that bound the sum
(\ref{sum}) in that case.
\end{itemize}
Lemma \ref{imp} directly follow from (\ref{id}). \qed
\end{proof}

\noindent {\bf Remarks.} We have in fact proved the following
sharper result. There exists a universally constant $C$ such that
\begin{align}
  \left(\frac{n}{{m_n}}\right)^{q/2}\|F\!L_n(x)\|_q^q \leq
  \left\{\begin{array} {ll}  (Cq)^q    & \mbox{in the Riemaniann case,}\\
              (Cq^{1+1/b}\sqrt n)^q    & \mbox{in the geometric case.}
\end{array}\right.
 \label{res}
\end{align}

\begin{proof}[Proof of Lemma \ref{liii}]
  The cases of Riemannian or geometric decay of the dependence
  coefficients are considered separately.
\begin{itemize}
\item{\bf Geometric decay} \ We present a technical lemma useful to
  deduce exponential probabilities from moment inequalities at any
  even order.
\begin{lemma}\label{momexp}
  If the variables $\{V_n\}_{n\in\mathbb{Z}}$ satisfies, for all
  $k\in\mathbb{N}^*$
\begin{equation}\label{hyp2}\|V_n\|_{2k}\le \phi(2k)\;,\end{equation}
where $\phi$ is an increasing function with $\phi(0)=0$.
Then:$$\displaystyle \mathbb{P}(|V_n|\ge\epsilon)\le
e^2\exp\left(-\phi^{-1}(\epsilon/e)\right)\;.$$
\end{lemma}
\begin{proof}
  By Markov's inequality and Assumption (\ref{hyp2}), we obtain
$$
\mathbb{P}\left(|V_n|\ge\epsilon\right)
\le\left(\frac{\phi(2k)}{\epsilon}\right)^{2k}\;.
$$
With the
convention $0^0=1$, the inequality is true for all $k\in\mathbb{N}$.
Reminding that $\phi(0)=0$, there exists an integer $k_0$ such that
$\phi(2k_0)\le \epsilon/e <\phi(2(k_0+1))$. Noting $\phi^{-1}$ the
generalized inverse of $\phi$, we have:
\begin{eqnarray*}
\mathbb{P}\left(|V_n|\ge\epsilon\right)&\le&\left(\frac{\phi(2k_0)}{\epsilon}\right)^{2k_0}
\le e^{-2k_0}=e^2e^{-2(k_0+1)}\;,\\
&\le&e^2\exp\left(-\phi^{-1}(\epsilon/e)\right)\;.
\end{eqnarray*}
\qed
\end{proof}

\noindent We rewrite the inequality (\ref{res}):
$\left\|\sqrt{\frac{n}{m_n}}F\!L_n\right\|_{2k}\le\phi(2k)$ with
$\displaystyle\phi(x)= Cx^\frac{b+1}{b}$ for a convenient constant
$C$. Applying Lemma \ref{momexp} to
$V_n=\sqrt{\frac{n}{m_n}}F\!L_n$ we obtain:
$$\mathbb{P}\left(|F\!L_n|\ge\epsilon\sqrt{\frac{{m_n}}n}\right)\le
e^2\exp\left(-\phi^{-1}(\epsilon/e)\right)\;,$$ and we obtain the
result of the Lemma \ref{liii}.

\item{\bf Riemannian decay}
In this case, the result of Lemma \ref{imp} is obtained only for
some values of $q$ depending of the value of the parameter $a$:
\begin{itemize}
\item In the case of $\eta$-dependence:
$$a>\max\left(q-1,\frac{1+\delta+2/d}{1-\delta},2+\frac{1}{d}\right)\;.$$
\item In the case of $\tilde \phi$-dependence:
$$a>\max\left(q-1,1+\frac{2}{d(1-\delta)},1+\frac{1}{d}\right)\;.$$
\end{itemize}

\noindent We consider that the assumptions of the Lemma \ref{liii}
on $a$ are satisfied in both cases of dependence. Then
$q_0=2\left\lceil\frac{a-1}{2}\right\rceil$ is the even integer such
that $a-1\le q_0 <  a+1$. It is the largest order such that the
assumptions of Lemma \ref{imp} (recalled above) are verified and
then the Lemma \ref{imp} gives us directly the rate of the moment:
$\displaystyle\lim_{n\to\infty}
\sup\left(\frac{n}{m_n}\right)^{q_0/2}\|
F\!L_n(x)\|^{q_0}_{q_0}<+\infty.$ We apply Markov to obtain the
result of Lemma \ref{liii}:
$$
\mathbb{P}\left(|F\!L_n(x)|\ge
  \epsilon\sqrt{\frac{{m_n}}{n}}\right) \le
\frac{\left(\sqrt{\frac{n}{{m_n}}}\|F\!L_n(x)\|_{q_0}\right)^{q_0}}{\epsilon^{q_0}}\;.
$$
\end{itemize}
\qed
\end{proof}

\begin{proof}[Proof of Lemma~\ref{rat}]
  We follow here Liebscher's strategy as in \cite{ang2}. We recover
  $B:= B(0,M)$, the ball of center $0$ and radius $M$, by at least
  $(4M\mu+1)^d$ balls $B_j=B(x_j,1/\mu)$. Then, under the assumption
  that $K_m(.,y)$ is supported on a compact of diameter proportional
  smaller than $1/m^{1/d}$, we have, for all $j$:
\begin{eqnarray} \label{sup}
\sup_{x\in B_j}|F\!L_n(x)|&\le&|\hat f_n(x_j)-\mathbb{E} \hat
f_n(x_j)|+C\frac{{m_n}^{1/d}}{\mu}(|\tilde f_n(x_j)-\mathbb{E}
\tilde f_n(x_j)|\\
\nonumber &&+2|\mathbb{E} \tilde f_n(x_j)|)\;,
\end{eqnarray}
with $C$ a constant and $\tilde f_n(x)=\frac
{1}{n}\sum_{i=1}^n\tilde K_{m_n}(x,X_i)$ where $\tilde K_{m_n}$ is
a kernel of type $\tilde
K_{m}(x,y)=K_0mk(x,y,x_j,1/m^{1/d})$. The $1/b$-Lipschitz function $k(x,y,a,b)$ is equal to $1$ on $B(a,b)$ and null outside $B(a,b+1/b)$. The constant $K_0$ is fixed in order that $\tilde K_{m_n}$ satisfies
properties (a), (b) and (c) of section \ref{nes}. Then using
(\ref{sup}) and with obvious short notation:
\begin{eqnarray*}
\mathbb{P}\left( \sup_{\|x\|\le
M}|F\!L_n(x)|>\epsilon\sqrt\frac{{m_n}}{n}\right)&\le&\sum_{j=1}^{(4M\mu+1)^d}\mathbb{P}\left(
\sup_{x\in B_j}|F\!L_n(x)|>\epsilon\sqrt\frac{{m_n}}{n}\right)\;,\\
&\le&(4M\mu+1)^d\left[\sup_{x\in B_j}\mathbb{P}\left(|F\!L_n(x_j)|>\epsilon\sqrt{\frac{{m_n}}{n}}\right)\right.\\
&&\left.+
\mathbb{P}\left(C\frac{m_n^{1/d}}{\mu}|\tilde{F\!L}_n(x)|>\epsilon\sqrt\frac{m_n}{n}\right)\right.\\
&&\left.+ \mathbb{P}\left(2C\frac{m_n^{1/d}}{\mu}|\mathbb{E} \tilde
f_n(x_j)|>\epsilon\sqrt{\frac{m_n}{n}}\right)\right]\;.
\end{eqnarray*}
Using the fact that $f$ is bounded, $\mathbb{E} \tilde f_n=\int \tilde
K_{m_n}(x_j,s) f(s)ds$ is bounded independently of $n$. Then taking $\mu= m_n^{1/d-1/2}n^{1/2}L(n)/\epsilon$ ensures that $\mathbb{P}\left(2C\frac{m_n^{1/d}}{\mu}|\mathbb{E} \tilde
f_n(x_j)|>\epsilon\sqrt{\frac{m_n}{n}}\right)$ is null for $n$ sufficiently large.  Applying Lemma
\ref{liii} on $f$ and $\tilde f$, uniform probability
inequality in both cases of geometric and Riemannian decays become:
\begin{eqnarray}\label{pi2}\mathbb{P}\left(\sup_{\|x\|\le
M}|F\!L_n(x)|\ge \epsilon_n\sqrt{\frac{m_n}{n}}\right)&\preceq&
\mu^d\exp\left(\
-C\epsilon_n^{\frac{b}{b+1}}\right)\;,\\
\label{pi1}\mathbb{P}\left(\sup_{\|x\|\le
M}|F\!L_n(x)|\ge \epsilon_n\sqrt{\frac{m_n}{n}}\right)&\preceq&
\mu^d\epsilon_n^{-q_0}\;.\end{eqnarray}

In the geometric case, fix $\epsilon_n$ as
$G(\log n)^{(b+1)/b}$ such that the bound becomes $\mu^dn^{-GC}$. Reminding that $\mu\le n$, the
sequence $\mu^d n^{-GC}$, bounded by $n^{d-GC}$, is summable for a
conveniently chosen constant $G$. Borel-Cantelli's Lemma then
concludes the proof in this case.

\noindent In the Riemannian case, take $\displaystyle\epsilon_n=
(m_n^{1-d/2}n^{1+d/2})^{\frac{1}{q_0+d}}\log n$ such that the bound becomes $\displaystyle n^{-1}\log^{-q_0}nL(n)$. Reminding that
$q_0\ge 2$, this sequence is summable and here again we conclude by
applying Borel-Cantelli's Lemma.  \qed
\end{proof}

\subsubsection{Proof of the theorems}

The order of magnitude of the bias is given by Assumption (d) and
the Lemmas provide bounds for fluctuation term. There only remain to
determine the optimal bandwidth $m_n$ in each case.

\begin{proof}[Proof of Theorem~\ref{theorem:1}]
  Applying Lemma \ref{imp} yields Theorem~\ref{theorem:1} when $q$ is
  an even integer. For any real $q$, Lemma \ref{imp} with $2(\lceil
  q/2 \rceil+1)\ge2$ and Jensen's inequalities yields:
\begin{eqnarray*}
\left(\frac{n}{m_n}\right)^{q/2}\mathbb{E}
|F\!L_n(x)|^q&=&\left(\frac{n}{m_n}\right)^{q/2}\mathbb{E} \left(
F\!L_n(x)^{2(\lceil
q/2 \rceil+1)}\right)^{q/\{2(\lceil q/2 \rceil+1)\}}  \;,
\\
& \le & \left(\left(\frac{n}{m_n}\right)^{\lceil q/2 \rceil+1}\mathbb{E}
F\!L_n(x)^{2(\lceil q/2 \rceil+1)}\right)^{q/\{2(\lceil q/2\rceil+1)\}} \;.
\end{eqnarray*}
Plugging this bound and the bound for the bias in~(\ref{fb}),
we obtain a bound for the $\mathbb{L}^q$-error of estimation:
\begin{equation*}
\|\hat f_n(x)-f(x)\|_q\le\|F\!L_n(x)\|_q+|R_n(x)|={O}
\left(\sqrt{\frac{m_n}{n}}+m_n^{-\rho/d}\right)\;.
\end{equation*}
The optimal bandwidth $ m^*_n=n^{\frac d{2\rho+d}}$ is the same as in the
i.i.d. case.  Thus [H6] holds with $\delta=\frac d{2\rho+d}$. For this
valued of $\delta$, the conditions on the parameter $a$ of
Lemma~\ref{liii} are equivalent to those of Theorem~\ref{theorem:1}.
\qed
\end{proof}

\begin{proof}[Proof of Theorem~\ref{thorem:2}]
  Applying the probability inequality (\ref{pi2}) in the proof of
  Lemma~\ref{rat} and the identity
  $\mathbb{E}|Y|^q=\int_0^{+\infty}\mathbb{P}\left(|Y|\ge
    t^{1/q}\right)dt$, we obtain
  $$
  \mathbb{E}\sup_{\|x\|\le M}| \hat f_n(x)-f(x)|^q = {O} \left(
    \left\{ \sqrt{\frac{m_n}{n}} \log^{{(b+1)/b}} (n) \right\}^q +
    m_n^{-q\rho/d}\right)\;.
  $$
  Lemma~\ref{rat} gives the rate of almost sure convergence:
  $$
  \sup_{\|x\|\le M}|\hat f_n(x)-f(x)|=_{a.s.}{O}
  \left(\sqrt{\frac{m_n}{n}}\log^{\frac{b+1}{b}}
    n+m_n^{-\rho/d}\right)\;.
  $$
  In both cases, the optimal bandwidth is $m^*_n = (n /
  \log^{2(b+1)/b}(n) )^{d/(2\rho+d)}$, which yields the rates claimed
  in Theorem~\ref{thorem:2}. \qed
\end{proof}

\begin{proof}[Proof of Theorem~\ref{theorem:3}]
  Applying the probability inequality (\ref{pi1}) and the same line of
  reasoning as in the previous proof, we obtain
\begin{equation*}
\mathbb{E}\sup_{\|x\|\le M}|\hat f_n(x)-f(x)|^q={O}
\left(\left(\sqrt{\frac{m_n}{n}}m_n^{\frac2{q_0+d}}\right)^q+m_n^{-q\rho/d}\right)\;,
\end{equation*}
where $q_0=2\left\lceil (a-1)/2 \right\rceil$. The optimal bandwidth $
m^*_n=n ^{d/(d+2\rho+2d/(q_0+d))}$ implies [H6] with
$\delta=d/(d+2\rho+2d/(q_0+d))$. For this value of $\delta$, the
conditions on $a$ of Lemma~\ref{liii} are satisfied as soon as $a\ge4$
and $\rho> 2d$.

\noindent  Lemma \ref{rat} gives the rate for the fluctuation
in the almost sure case. This leads the optimal bandwidth
$$
m^*_n=\left(n^{q_0-2}/\log^{q_0+d}(n) \right)^{\frac{d}{d(q_0+2)+\rho(q_0+d)}}.
$$
We then deduce the two different rates of Theorem~\ref{theorem:3},
either in the almost sure or in the $\mathbb{L}^q$ framework.  \qed
\end{proof}


\end{document}